\documentclass[11pt,twoside]{article}
\usepackage[mathscr]{euscript}
\usepackage{secdot}
\usepackage[utf8]{inputenc}
\usepackage{amsmath, amsthm, amscd, amsfonts, amssymb, color}
\usepackage[bookmarksnumbered,plainpages]{hyperref}
\allowdisplaybreaks

\date{}

\begin{document}

\centerline{}

\centerline {\Large{\bf Fixed Point Theorems for TSR-Contraction Mapping}}
\centerline{}
\centerline {\Large{\bf in Probabilistic Metric Spaces}}
\centerline{}
\centerline{\textbf{Sanjay Roy}}
\centerline{Department of Mathematics, Uluberia College,}
\centerline{Uluberia, Howrah,  West Bengal, India}
\centerline{E-mail: sanjaypuremath@gmail.com}
\centerline{}

\centerline{\textbf{T. K. Samanta}}
\centerline{Department of Mathematics, Uluberia College,}
\centerline{Uluberia, Howrah,  West Bengal, India}
\centerline{Email:~mumpu$_{-}$tapas5@yahoo.co.in}
\newcommand{\mvec}[1]{\mbox{\bfseries\itshape #1}}
\centerline{}
\newtheorem{Theorem}{\quad Theorem}[section]

\newtheorem{definition}[Theorem]{\quad Definition}

\newtheorem{theorem}[Theorem]{\quad Theorem}

\newtheorem{remark}[Theorem]{\quad Remark}

\newtheorem{corollary}[Theorem]{\quad Corollary}

\newtheorem{note}[Theorem]{\quad Note}

\newtheorem{lemma}[Theorem]{\quad Lemma}

\newtheorem{example}[Theorem]{\quad Example}

\newtheorem{notation}[Theorem]{\quad Notation}

\newtheorem{result}[Theorem]{\quad Result}
\newtheorem{conclusion}[Theorem]{\quad Conclusion}

\newtheorem{proposition}[Theorem]{\quad Proposition}
\newtheorem{prop}[Theorem]{\quad Property}

\begin{abstract}
\textbf{\emph{The concept of fixed point plays a crucial role in various fields of applied mathematics. The aim of this paper is to establish the existence of a unique fixed point of some type of functions which satisfy a new contraction principle, namely, TSR-contraction principle in various types of probabilistic metric spaces. The proposed contraction mapping is different from our traditional definitions of contraction mapping. 
}}
\end{abstract}

\section{Introduction}

  The distance between two points $x$ and $y$ in a probabilistic metric space is a distance distribution function $F_{x,\, y}$ whose value at a real number $t$ is interpreted as the probability that the distance between the points $x$ and $y$ is less than equal to $t$. In 1942, K. Menger \cite{Menger} played important role to develop the concept of probabilistic metric space. Then in 1958, B. Schweizer, A. Sklar \cite{Schweizer1, Schweizer2, Schweizer3} defined this space. With the help of this definition, many works \cite{Hadzic13, Yaari, Zhou} have been done on fixed point theorems for various types of contraction mappings.
 
 In 1983, S. S. Chang \cite{Chang1} introduced a contraction mapping with the help of a non-decreasing function $\Phi$ and found out the existence of common random fixed point in which the two sequences of $d$-continuous random operators converge. Then he proved the common fixed point theorems  for multi-valued mappings in Menger PM-spaces \cite{Chang2}.
 
 Since 1980, many works had been done by O. Hadzic to develop this concept of fixed point theorem as follows: 
In 1980, O. Hadzic \cite{Hadzic5} introduced almost fixed point theorem for multivalued mappings in a topological vector space.
In 1982, O. Hadzic \cite{Hadzic6} developed the concept of fixed point theorems for multivalued functions in a topological vector space.
In 1982, O. Hadzic \cite{Hadzic8} presented condition of existence of a unique fixed point of sum of two functions.
In 1988, O. Hadzic \cite{Hadzic11} investigated the common fixed point of two functions in a probabilistic metric space with convex structure.
In 1990, O. Hadzic \cite{Hadzic12} gave sufficient conditions for existence of fixed point of a multivalued B-contraction mapping.

In 1983, T. L. Hicks \cite{Hicks1} introduced a contraction mapping in probabilistic metric space $(X, \mathcal{F}, \tau)$ as follows: for $t>0$, $F_{f(x),\;f(y)}(kt)\geq 1-kt$ whenever $F_{x, y}(t)>1-t$, where $x,\, y\in X$ and developed a few results on fixed point theorem.
In 2007, S. Kutukcu etl.\cite{Kutukcu} found out the common fixed point for self maps in Menger probabilistic metric spaces.
In 2013, A. Mbarki etl.\cite{Mbarki} found out common fixed point in probabilistic metric spaces for four mappings.
In 2013, R. Singh etl.\cite{Singh} found out the common fixed point for four mappings with weak compatibility in probabilistic metric spaces. 
In 2014, N. Sharma \cite{Sharma} proved the existence of fixed points  for a self-map on a Menger space and found out the unique fixed point for a sequence of self maps on a complete probabilistic metric space.
 In 2014, C. Zhou etl. \cite{Zhou} generalized the concept of the Menger probabilistic metric space and investigated some topological properties of this space.
In 2022, A. A. Yaari etl. \cite{Yaari} showed a fixed point for this continuity map between a complete Hausdorff cone metric space and quasi-cone 
metric space.
  
  In this paper, we have introduced another type of  contraction mapping and established fixed point theorem on Menger space. Here we have considered probabilistic metric space with a particular type of triangle function, as a result we have got another new set of sufficient conditions for establishing fixed point theorem in probabilistic metric space.
  
\section{Preliminaries}
\smallskip\hspace{.6 cm} In this section, we recall some definitions and results in Probabilistic metric space.
\begin{definition}\cite{Schweizer3}
A binary relation $\ast$ on $[0,1]$ is said to be a triangular norm or $t-$norm if the following condition are satisfied:\\
$(i)$  $x\ast y=y\ast x$,\\
$(ii)$  $x\ast(y\ast z)=(x\ast y)\ast z$,\\
$(iii)$  $x\ast y\leq x\ast z$ whenever $y\leq z$\\
$(iv)$  $x\ast 1= x$ for all $x, y, z\in [0, 1]$.
\end{definition}

\begin{definition}
A distribution function is a function $F:[-\infty,\, \infty]\rightarrow [0, 1]$ which is left continuous on $\mathbb{R}$, non-decreasing and $F(-\infty)=0$, $F(\infty)=1$.
\end{definition}

\begin{definition}
The Dirac distribution function $H_a:[-\infty,\,\infty]\rightarrow [0, 1]$ is defined for $a\in[-\infty,\, \infty)$ by\\
\begin{center}
$H_a(u) = 
\begin{cases} 
0 & \text{if } u \in [-\infty,\, a], \\
1 & \text{if } u\in(a,\,\infty], 
\end{cases}$
\end{center}
and for $a=\infty$ by\\
\begin{center}
$H_\infty(u) = 
\begin{cases} 
0 & \text{if } u \in [-\infty,\, \infty), \\
1 & \text{if } u=\infty. 
\end{cases}$
\end{center}
\end{definition}

\begin{definition}
A distance distribution function $F:[-\infty,\, \infty]\rightarrow [0, 1]$ is a distribution function with $f(0)=0$. The family of all distance distribution functions is denoted by $\bigtriangleup^+$. 
\end{definition}

\begin{definition}
A triangle function $\tau$ is a binary operation on $\bigtriangleup^+$ that is commutative, associative and non-decreasing in each place and has the identity $H_0$.
\end{definition}

\begin{example}
Let $\ast$ be a left-continuous t-norm. Then the function\\  $\textbf{T}: \bigtriangleup^+ \times \bigtriangleup^+ \rightarrow \bigtriangleup^+$ defined by $\textbf{T}(F, G)(x)= F(x)\ast G(x) $ is a triangle function.
\end{example}

\begin{definition}\cite{Hadzic14}
A probabilistic metric space (briefly, a PM space)  is a triple $(X, \mathcal{F},\tau)$ where $X$ is a non empty set and $\mathcal{F}: X\times X\rightarrow \bigtriangleup^+$, given by $(p, q)\mapsto F_{p, q}$ and $\tau$ is a triangle function such that the following conditions are satisfied for all $p,q,r$ in $X$:\\
$(i)$ $F_{p, q}= H_0$ if and only if $p=q$,\\
$(ii)$ $F_{p, q}=F_{q, p}$,\\
$(iii)$ $F_{p, r}\geq \tau\left( F_{p, q},\;F_{q,\,r}\right) $.\\
\end{definition}

\begin{definition}\label{D2}\cite{Hadzic14}
Let $(X, \mathcal{F}, \tau)$ be a probabilistic metric space and $\tau=\tau_\ast$, where
$\tau_\ast\left( F,\;G\right)(t)=\sup\{ F(t_1)\ast G(t_2):\;t_1+t_2=t, \; t_1, t_2\geq 0\}$ for all $t\geq 0$ for a t-norm $\ast$. Then $(X, \mathcal{F}, \tau)$ is called a Menger Space which will be denoted by $(X, \mathcal{F}, \tau_\ast)$.
\end{definition}

\begin{remark}\cite{Hadzic14}
If the t-norm $\ast$ is left continuous then $\tau_\ast$ of the Definition \ref{D2} is a triangle function. Then we have
$F_{p,\, r}(t_1+t_2)\geq \tau_\ast\left( F_{p, q},\;F_{q, r}\right)(t_1+t_2)\geq F_{p,\, q}(t_1)\ast F_{q,\,r}(t_2)$
 for all $t_1,\, t_2\geq 0$ and for all $p, q, r\in X$. 
 \end{remark}

\begin{definition}\cite{Hicks1}
A mapping $f$ of a PM- space  $(X, \mathcal{F}, \tau)$ into itself is a contraction mapping if there exists $k\in (0, 1)$ such that for each $x, y\in X$, $F_{f(x),\;f(y)}(kt)\geq F_{x, y}(t)$ for all $t> 0$. If $f$ satisfies the above Contraction principle then $f$ is called $B$-contraction mapping.

Another generalization of the contraction mapping in PM -space is that for $t>0$, $F_{f(x),\;f(y)}(kt)\geq 1-kt$ whenever $F_{x, y}(t)>1-t$. If $f$ satisfies the above Contraction principle then $f$ is called $H$-contraction mapping.
\end{definition}

\begin{definition}\cite{Schweizer2}
Let $(X, \mathcal{F}, \tau)$ be a probabilistic metric space. A sequence $\{x_n\}$ in $X$ is said to converge to a point $x\in X$ if for any $\alpha\in (0, 1)$ and for any $t>0$, there exists a natural number $m=m(\alpha, t)$ such that $F_{x_{n},\, x}(t) > 1-\alpha$ for all $n\geq m$
\end{definition}

\begin{definition}
Let $(X, \mathcal{F}, \tau)$ be a probabilistic metric space. A sequence $\{x_n\}$ in $X$ is said to be a Cauchy sequence if for any $\alpha\in (0, 1)$ and for any $t>0$, there exists a natural number $m=m(\alpha, t)$ such that $F_{x_{n+p},\, x_n}(t) > 1-\alpha$ for all $n\geq m$.
\end{definition}

\begin{definition}
Let $(X, \mathcal{F}^1, \tau_1)$ and $(Y, \mathcal{F}^2, \tau_2)$ be two probabilistic metric spaces and $f: X\rightarrow Y$ be a mapping. Then $f$ is said to be continuous at a point $a\in X$ if for any $\epsilon\in (0, 1)$ and $t> 0$, there exists $\delta=\delta(\epsilon, a, t)\in (0, 1)$ such that $F^1_{x, a}(t)> 1-\delta$ implies $F^2_{f(x), f(a)}(t)>1-\epsilon$.
\end{definition}

\begin{theorem}\label{th21}
Let $(X, \mathcal{F}^1, \tau_1)$ and $(Y, \mathcal{F}^2, \tau_2)$ be two PM-spaces and $f: X\rightarrow Y$ be a continuous mapping. If  $x_n\rightarrow a$ as $n\rightarrow\infty$ in $(X, \mathcal{F}^1, \tau_1)$ then $f(x_n)\rightarrow f(a)$ as $n\rightarrow\infty$ in $(Y, \mathcal{F}^2, \tau_2)$.
\end{theorem}

\proof Let $\epsilon\in (0, 1)$ and $t> 0$. 
Since $f$ is continuous at $a$, there exists $\delta\in (0, 1)$ such that $F^1_{x, a}(t)> 1-\delta$ implies that $F^2_{f(x),\, f(a)}(t)>1-\epsilon$. Again,
 Since $x_n\rightarrow a$ as $n\rightarrow\infty$, there exists $m\in \mathbb{N}$ such that $F^1_{x_n,\, a}(t)>1-\delta$ for all $n\geq m$.
Therefore, for all $n\geq m$, we have $F^2_{f(x_n),\, f(a)}(t)>1-\epsilon$. Hence $f(x_n)\rightarrow f(a)$ as $n\rightarrow\infty$ in $(Y, \mathcal{F}^2, \tau_2)$.

\begin{theorem}
Let $(X, \mathcal{F}, \tau_{\ast})$ be a Menger space and   $\alpha\ast\alpha\geq\alpha$ for all $\alpha\in[0, 1]$. Then every convergent sequence of $(X, \mathcal{F}, \tau_{\ast})$ has a unique limit in $(X, \mathcal{F}, \tau_{\ast})$.
\end{theorem}

\proof Assume that the sequence $\{x_n\}$   converges to points $x$ and $y$ in $(X, \mathcal{F}, \tau_{\ast})$. Let $\alpha\in (0, 1)$ and $t> 0$. Then there exists $m\in \mathbb{N}$ such that $F_{x_n,\, x}(\frac{t}{2})>1-\alpha$ and  $F_{x_n,\, y}(\frac{t}{2})>1-\alpha$ for all $n\geq m$. \\
 Now $F_{x,\, y}(t)\geq \tau_{\ast}\left( F_{x,\, x_n},\;F_{x_n,\, y}\right)(t)\geq F_{x,\, x_n}(\frac{t}{2})\ast F_{x_n,\, y}(\frac{t}{2}) \geq (1-\alpha)\ast (1-\alpha)\geq (1-\alpha)$ for all $n\geq m$. So $F_{x,\, y}(t)= 1$ as $\alpha$ is chosen arbitrarily. Therefore  $F_{x,\, y}(t)= 1$ for all $t> 0$. Hence $F_{x,\, y}=H_0$, that is, $x=y$.

\section{Fixed point theorem on Menger Space}
\begin{definition}\label{D1}
Let $(X, \mathcal{F}, \tau)$ be a probabilistic metric space. A mapping $f:X\rightarrow X$ is called probabilistic TSR contraction mapping or, TSR contraction mapping if for $t> 0$, there exists $k\in (0, 1)$ such that $1-F_{f(x), f(y)} (kt)\leq k(1-F_{x, y}(t))$ for every $x, y\in X$. Here $k$ will be called contraction constant.

From the Definition \ref{D1}, it follows that for all $\alpha\in (0, 1)$ with $F_{x, y}(t)> \alpha$ implies $1-F_{f(x),\, f(y)}(kt)<k(1-\alpha)$, that is, $F_{f(x),\, f(y)}(kt)>1-k(1-\alpha)>\alpha$.
\end{definition}

\begin{definition}
Let $(X, \mathcal{F}, \tau)$ be a probabilistic metric space. A mapping $f:X\rightarrow X$ is called probabilistic TSR P-contraction mapping or, TSR P-contraction mapping if for $t> 0$, there exists $k\in (0, 1)$ such that $1-F_{f(x), f(y)} (k^{m+1}t)\leq k(1-F_{x, y}(k^mt))$ for every $x, y\in X$ and for all $m\in\mathbb{N}$.
\end{definition}

\begin{theorem}\label{th31}
Let $(X, \mathcal{F}, \tau_{\ast})$ be a complete Menger space and $f: X\rightarrow X$ be a probabilistic TSR P-contraction mapping with contraction constant less than equal to  $\frac{1}{2}$. If $\alpha\ast\alpha\geq\alpha$ for all $\alpha\in[0, 1]$ then $f$ has a unique fixed point in $X$.\end{theorem}

\proof Let $t> 0$. Since $f$ is a probabilistic TSR P-contraction mapping with contraction constant less than equal to  $\frac{1}{2}$, there exists $k\in (0, \frac{1}{2}]$ such that\\
 $1-F_{f(a), f(b)} (k^{r+1}t)\leq k(1-F_{a, b}(k^r t))$ for all $a, b\in X$ and for all $r\in\mathbb{N}$.\\
   Let $x_0\in X$ and $x_n=f^n(x_0)$ for all $n\in\mathbb{N}$. Then taking $a=x_1=f(x_0)$ and $b=x_0$ we get
\begin{align*}
&&1-F_{f(f(x_0)), f(x_0)} (kt)&\leq k(1-F_{f(x_0), x_0}(t)),\\ 
&\text{or}, &1-F_{x_2, x_1} (kt)&\leq k\left(1-F_{x_1, x_0}(t)\right) ,\\
& \text{or}, &F_{x_2, x_1} (kt)&\geq 1-k\left(1-F_{x_1, x_0}(t)\right).\\
&\text{or}, &-F_{x_2, x_1} (kt)&\leq -1+k\left(1-F_{x_1, x_0}(t)\right).\\
&\text{or}, &1-F_{x_2, x_1} (kt)&\leq k\left(1-F_{x_1, x_0}(t)\right).\\
&\text{or}, &k\left( 1-F_{x_2, x_1} (kt)\right) &\leq k^2\left(1-F_{x_1, x_0}(t)\right).\\
&\text{or}, &1-k\left( 1-F_{x_2, x_1} (kt)\right) &\geq 1-k^2\left(1-F_{x_1, x_0}(t)\right). \smallskip\hspace{2cm}\cdots (i)
\end{align*} 
 Again taking $a=x_2$ and $b=x_1$ we get
  \begin{align*}
  && 1-F_{x_3, x_2} (k^2t)&\leq k\left( 1-F_{x_2, x_1}(kt)\right),\\
  &\text{or}, &F_{x_3, x_2} (k^2t)&\geq 1-k\left( 1-F_{x_2, x_1}(kt)\right) \geq 1-k^2\left(1-F_{x_1, x_0}(t)\right).\\
&\text{or}, &-F_{x_3, x_2} (k^2t)&\leq -1+k^2\left(1-F_{x_1, x_0}(t)\right).\\
&\text{or}, &1-F_{x_3, x_2} (k^2t)&\leq k^2\left(1-F_{x_1, x_0}(t)\right).\\
&\text{or}, &k\left( 1-F_{x_3, x_2} (k^2t)\right) &\leq k^3\left(1-F_{x_1, x_0}(t)\right).\\
&\text{or}, &1-k\left( 1-F_{x_3, x_2} (k^2t)\right) &\geq 1-k^3\left(1-F_{x_1, x_0}(t)\right).
\end{align*}
So, by induction we get $F_{x_{n+1},\, x_n} (k^nt)\geq 1- k^n\left( 1-F_{x_1,\, x_0}(t)\right) $ for all $n\in\mathbb{N}$.\\
We now show that $\{x_n\}$ is a Cauchy sequence in $(X, \mathcal{F}, \tau_{\ast})$. Let $\alpha\in (0, 1)$. Then there exists a natural number $m$ such that $ k^{m}\left( 1-F_{x_1,\, x_0}(t)\right)<  \alpha$. Let $p $ be any natural number. Now,
\begin{align*}
 F_{x_{n+p}, \,x_n}(t)&=F_{x_n, \,x_{n+p}}(t)\\
&\geq \tau_{\ast}\left( F_{x_n, \,x_{n+1}},F_{x_{n+1}, \,x_{n+p}}\right)(t)\\
&\geq F_{x_{n}, \,x_{n+1}}(\frac{t}{2})\ast F_{x_{n+1}, \,x_{n+p}}(\frac{t}{2})\\
&\geq  F_{x_{n}, \,x_{n+1}}(\frac{t}{2})\ast \tau_{\ast}\left( F_{x_{n+1}, \,x_{n+2}}, F_{x_{n+2}, \,x_{n+p}}\right)(\frac{t}{2}),\; \text{as}\; \ast \;\text{is non-decreasing.}\\
&\geq F_{x_{n}, \,x_{n+1}}(\frac{t}{2})\ast  F_{x_{n+1}, \,x_{n+2}}(\frac{t}{2^2})\ast F_{x_{n+2}, \,x_{n+p}}(\frac{t}{2^2}), \;\text{as}\; \ast \;\text{is associative.}\\
&\vdots\\
&\geq F_{x_{n}, \,x_{n+1}}(\frac{t}{2})\ast  F_{x_{n+1}, \,x_{n+2}}(\frac{t}{2^2})\ast \cdots \ast F_{x_{n+p-1}, \,x_{n+p}}(\frac{t}{2^p})\\
&\geq F_{x_{n}, \,x_{n+1}}(k^{n}t)\ast  F_{x_{n+1}, \,x_{n+2}}(k^{n+1}t)\ast \cdots \ast F_{x_{n+p-1}, \,x_{n+p}}(k^{n+p-1}t)\; \text{as} \;k\leq\frac{1}{2}\\
&= F_{x_{n+1}, \,x_{n}}(k^{n}t)\ast  F_{x_{n+2}, \,x_{n+1}}(k^{n+1}t)\ast \cdots \ast F_{x_{n+p}, \,x_{n+p-1}}(k^{n+p-1}t)\\
&\geq \left( 1- k^{n}\left( 1-F_{x_1,\, x_0}(t)\right)\right)  \ast \left( 1- k^{n+1}\left( 1-F_{x_1,\, x_0}(t)\right) \right) \ast \cdots \\
& \smallskip\hspace{.5 cm}\ast \left( 1- k^{n+p-1}\left( 1-F_{x_1,\, x_0}(t)\right)\right),  \text{ as}\; k\in \left( 0, \frac{1}{2}\right]  \;\text{and}\; \ast\; \text{is non-decreasing.}\\
&\geq 1- k^{n}\left( 1-F_{x_1,\, x_0}(t)\right).
\end{align*}
So, $F_{x_{n+p}, \,x_n}(t)>1-\alpha$ whenever $1- k^{n}\left( 1-F_{x_1,\, x_0}(t)\right)> 1-\alpha$.\\
Now
$1- k^{n}\left( 1-F_{x_1,\, x_0}(t)\right)> 1-\alpha$ iff $ k^{n}\left( 1-F_{x_1,\, x_0}(t)\right)< \alpha$ if  $n\geq m$, as\\
 $ k^{m}\left( 1-F_{x_1,\, x_0}(t)\right)< \alpha$.
Therefore we see that $F_{x_{n+p}, \,x_n}(t)>1-\alpha$ whenever $n\geq m$.
Thus $\{x_n\}$ is a Cauchy sequence in $(X, \mathcal{F}, \tau_{\ast})$. Since the space is complete Menger space, there exists $y\in X$ such that $\lim x_n=y$, that is, $\lim f^n (x_0)=y$. We now show that $y$ is a fixed point of the function $f$, that is, $f(y)=y$.  We have
\begin{align*}
F_{f(y),\, y}(t)&\geq \tau_{\ast}\left( F_{f(y),\, x_n},\, F_{x_n, y}\right) (t)\\
 &\geq F_{f(y),\, x_n}(\frac{t}{2})\ast F_{x_n, y}(\frac{t}{2})\\
 &=F_{f(y),\, f^n(x_0)}(\frac{t}{2})\ast F_{x_n, y}(\frac{t}{2})\\
& = F_{f(y),\, f(x_{n-1})}(\frac{t}{2})\ast F_{x_n, y}(\frac{t}{2})\\
& \geq F_{f(y),\, f(x_{n-1})}(kt)\ast F_{x_n, y}(\frac{t}{2})\;\text{as}\; k\leq\frac{1}{2}\\
 & \geq \left[ 1-k\left( 1-F_{y,\, x_{n-1}}(t)\right) \right] \ast F_{x_n, y}(\frac{t}{2}).\smallskip\hspace{2cm}\cdots (ii)
 \end{align*}
 Now since $\lim x_n = y$, for any $\alpha \in (0, 1)$, there exists $q\in \mathbb{N}$ such that $F_{y,\, x_{n-1}}(t)> 1-\alpha$ and $F_{x_n, y}(\frac{t}{2})> 1-\alpha$ for all $n\geq q$. So, from $(ii)$ we get \\
 $F_{f(y),\, y}(t)\geq ( 1-k\alpha) \ast (1-\alpha)$, as $\ast$ is non-decreasing.\\
 So, $F_{f(y),\, y}(t)\geq 1-\alpha$ for all $t>0$ and for all $\alpha\in (0, 1)$. Therefore $F_{f(y),\, y}(t)=1$ for all $t>0$. So $F_{f(y),\, y}= H_0$.
Thus $f(y)=y$. We now verify that there exists only one fixed point of $f$. Let $z\in X$ be such that $f(z)=z$. Then \\$F_{y, z}(t)\geq F_{y, z}(kt)= F_{f(y), f(z)}(kt)\geq 1-k\left( 1-F_{y, z} (t)\right)$ implies\\
 $1- F_{y, z}(t)\leq k\left( 1-F_{y, z} (t)\right)$. If  $F_{y, z}(t)< 1$, that is, if $1- F_{y, z}(t)> 0$, we obtain $k\geq 1$ which is not true. Thus $F_{y, z}(t)=1$ for all $t> 0$. So, $F_{y, z}= H_0$, that is, $y=z$. This completes the proof of the theorem.

\begin{definition}
Let $(X, \mathcal{F}, \tau)$ be a probabilistic metric space, $x_0\in X$, $r\in (0, 1)$ and $t>0$. The t-open sphere centered at $x_0$ of radius $r$  is denoted by $S_t(x_0,\, r,\,\mathcal{F} )$ and is defined by $S_t(x_0,\, r,\,\mathcal{F} )=\{x\in X: F_{x, \, x_0}(t)>1-r\}$. 
\end{definition}

\begin{definition}
Let $(X, \mathcal{F}, \tau)$ be a probabilistic metric space, $t>0$ and $S\subseteq X$. A point $x_0\in X$ is said to be a t-limit point of $S$ if for any $r\in (0, 1)$, the t-open sphere $S_t(x_0,\, r,\,\mathcal{F} )$  centered at $x_0$ of radius $r$ intersects $S$.
\end{definition}

\begin{definition}
Let $(X, \mathcal{F}, \tau)$ be a probabilistic metric space and $t>0$. A subset $S\subseteq X$ is said to be a t-closed in $X$ if $S$ contains all its t-limit points.
\end{definition}

\begin{definition}
Let $(X, \mathcal{F}, \tau)$ be a probabilistic metric space, $x_0\in X$, $r\in (0, 1)$ and $t> 0$. The t-closed sphere centered at $x_0$ of radius $r$ is denoted by $S_t(x_0,\, r,\,\mathcal{F} )$ and is defined by $\overline{S_t}(x_0,\, r,\,\mathcal{F} )=\{x\in X: F_{x, \, x_0}(t)\geq 1-r\}$. 
\end{definition}

\begin{theorem}
Let $(X, \mathcal{F}, \tau_{\ast})$ be a complete Menger space and $f: X\rightarrow X$ be a probabilistic TSR P-contraction mapping on $\overline{S_t}(x_0,\, r,\,\mathcal{F} )$ with contraction constant less than equal to $\frac{1}{2}$. If $\overline{S_t}(x_0,\, r,\,\mathcal{F} )$ is t-closed, $\alpha\ast\alpha\geq\alpha$ for all $\alpha\in[0, 1]$  and $F_{x_0,\, f(x_0)}(u)>1-r$ for all $u>0$ then $f$ has a unique fixed point in $\overline{S_t}(x_0,\, r,\,\mathcal{F} )$.
\end{theorem}

\proof Since $f$ be a probabilistic TSR P-contraction mapping on $\overline{S_t}(x_0,\, r,\,\mathcal{F} )$ with contraction constant less than equal to $\frac{1}{2}$, there exists $k\in (0, \frac{1}{2}]$ such that $1-F_{f(a), f(b)} (k^{m+1}\frac{t}{2})\leq k(1-F_{a, b}(k^m \frac{t}{2}))$ for all $a, b\in \overline{S_t}(x_0,\, r,\,\mathcal{F} )$ and for all $m\in\mathbb{N}$.\\
Let $x_1=f(x_0)$, $x_2=f(x_1)=f^2(x_0)$, $\cdots$,  $x_n=f(x_{n-1})=f^n(x_0)$.\\
Now $F_{x_1,\, x_0}(t)=F_{f(x_0),\, x_0}(t)>1-r$ by taking $u=t$. So, $x_1\in \overline{S}(x_0,\, r,\,\mathcal{F} )$. Assume that $x_1, x_2, \cdots, x_{n-1}\in \overline{S}(x_0,\, r,\,\mathcal{F} )$. We now show that $x_n\in \overline{S}(x_0,\, r,\,\mathcal{F} )$. We have
\begin{align*}
&&1-F_{x_2,\, x_1}(k\frac{t}{2})&=1-F_{f(x_1),\, f(x_0)}(k\frac{t}{2})\leq k\left( 1-F_{x_1,\, x_0}(\frac{t}{2})\right) ,\\
&\text{or}, &1- \frac{1-F_{x_2,\, x_1}(k\frac{t}{2})}{k}&\geq F_{x_1,\, x_0}(\frac{t}{2})>1-r,  \text{as}  F_{x_0,\, f(x_0)}(u)>1-r \text{ for all } u> 0,\\
&\text{or}, &\frac{1-F_{x_2,\, x_1}(k\frac{t}{2})}{k}&<r,\\
&\text{or}, &F_{x_2,\, x_1}(k\frac{t}{2})&>1-kr,\smallskip\hspace{2cm}\cdots (i)
\end{align*}
Again, 
\begin{align*}
&&1-F_{x_3,\, x_2}(k^2\frac{t}{2})&=1-F_{f(x_2),\, f(x_1)}(k^2\frac{t}{2})\leq k\left( 1-F_{x_2,\, x_1}(k\frac{t}{2})\right),\\
&\text{or}, &1- \frac{1-F_{x_3,\, x_2}(k^2\frac{t}{2})}{k}&\geq F_{x_2,\, x_1}(k\frac{t}{2})>1-kr,  \text{ by } (i),\\
&\text{or}, &\frac{1-F_{x_3,\, x_2}(k^2\frac{t}{2})}{k}&<kr,\\
&\text{or}, &F_{x_3,\, x_2}(k^2\frac{t}{2})&>1-k^2r.\\
\end{align*}
So, by induction we get $F_{x_{n+1},\, x_{n}}(k^n\frac{t}{2})>1-k^{n}r$ for all $n\in \mathbb{N}$.
Now,
\begin{align*}
F_{x_{n}, \,x_0}(t)&=F_{x_{0}, \,x_n}(t)\\
&\geq \tau_{\ast}\left( F_{x_{0}, \,x_{1}},F_{x_{1}, \,x_{n}}\right)(t)\\
&\geq F_{x_{0}, \,x_{1}}(\frac{t}{2})\ast F_{x_{n-1}, \,x_{0}}(\frac{t}{2})\\
&\geq  F_{x_{0}, \,x_{1}}(\frac{t}{2})\ast \tau_{\ast}\left( F_{x_{1}, \,x_{2}}, F_{x_{2}, \,x_{n}}\right)(\frac{t}{2}), \text{ as } \ast \text{ is non-decreasing.}\\
&\geq F_{x_{0}, \,x_{1}}(\frac{t}{2})\ast  F_{x_{1}, \,x_{2}}(\frac{t}{2^2})\ast F_{x_{2}, \,x_{n}}(\frac{t}{2^2}), \text{ as } \ast \text{ is associative.}\\
&\vdots\\
&\geq F_{x_{0}, \,x_{1}}(\frac{t}{2})\ast  F_{x_{1}, \,x_{2}}(\frac{t}{2^2})\ast \cdots \ast F_{x_{n-1}, \,x_{n}}(\frac{t}{2^n})\\
&= F_{x_{1}, \,x_{0}}(\frac{t}{2})\ast  F_{x_{2}, \,x_{1}}(\frac{t}{2^2})\ast \cdots \ast F_{x_{n}, \,x_{n-1}}(\frac{t}{2^n}) \\
&\geq F_{x_{1}, \,x_{0}}(\frac{t}{2})\ast  F_{x_{2}, \,x_{1}}(k\frac{t}{2})\ast \cdots \ast F_{x_{n}, \,x_{n-1}}(k^{n-1}\frac{t}{2}), \text{ as } k\leq\frac{1}{2}\\
&\geq \left( 1- r\right)   \ast \left( 1- kr \right) \ast \cdots \ast \left( 1- k^{n-1}r\right)\\
&\geq \left( 1- r\right)   \ast \left( 1- r \right) \ast \cdots \ast \left( 1- r\right), \text{ as } \ast \text{ is non-decreasing.}\\
&\geq 1- r.
\end{align*} 
Therefore, $x_n\in \overline{S_t}(x_0,\, r,\,\mathcal{F} )$. Thus, $x_n\in \overline{S}(x_0,\, r,\,\mathcal{F} )$ for all $n\in \mathbb{N}$. 
Let $u>0$. Since $f$ is a probabilistic TSR P-contraction mapping on $\overline{S_t}(x_0,\, r,\,\mathcal{F} )$, there exists $p\in (0, \frac{1}{2}]$ such that $1-F_{f(x),\, f(y)}(p^{m+1}u)\leq p\left( 1-F_{x,\, y}(p^mu)\right)$ for all $x,\, y\in \overline{S_t}(x_0,\, r,\,\mathcal{F} )$ and for all $m\in\mathbb{N}$.
Then by the above procedure, we get  $F_{x_{n+1},\, x_{n}}(p^nu)>1-p^{n}r$ for all $n\in \mathbb{N}.$
Now by the procedure of the proof of the Theorem \ref{th31}, we can show that the sequence $\{x_n\}$ is a Cauchy sequence. Since $(X, \mathcal{F}, \tau_{\ast})$ is a complete Menger space,  $\{x_n\}$ converges to an element $y$ of $X$. So, $y$ is a t-limit point of $\overline{S_t}(x_0,\, r,\,\mathcal{F} )$. Again Since  $\overline{S_t}(x_0,\, r,\,\mathcal{F} )$ is t-closed, $y\in \overline{S_t}(x_0,\, r,\,\mathcal{F} )$.  We now show that $y$ is a unique fixed point of $f$ in $\overline{S_t}(x_0,\, r,\,\mathcal{F} )$. Let $z$ be a fixed point of $f$ in $\overline{S_t}(x_0,\, r,\,\mathcal{F} )$. Then $F_{y, z}(u)\geq F_{y, z}(pu)= F_{f(y), f(z)}(pu)\geq 1-p\left( 1-F_{y, z} (u)\right)$ implies $1- F_{y, z}(u)\leq p\left( 1-F_{y, z} (u)\right)$. If  $F_{y, z}(u)< 1$, that is, if $1- F_{y, z}(u)> 0$, we obtain $p\geq 1$ which is not true. Thus $F_{y, z}(u)=1$. Since $u$ is chosen arbitrarily, $F_{y, z}= H_0$, that is, $y=z$. This completes the proof of the theorem. 

\section{Fixed point theorem on Probabilistic Metric Space}
Throughout this section we consider the probabilistic metric $($ PM $)$ space $(X, \mathcal{F}, \tau_\diamond)$ where $\diamond$ is a left-continuous t-norm and  $\tau_\diamond$ is defined as follows: $\tau_\diamond (F, G)(x)=F(x)\diamond G(x)$ for all $F,\, G\in\bigtriangleup^+$ and $x\geq 0$.

\begin{theorem}
Let $\alpha\diamond\alpha\geq\alpha$ for all $\alpha\in[0, 1]$. Then every convergent sequence of the PM Space $(X, \mathcal{F}, \tau_\diamond)$ has a unique limit in $(X, \mathcal{F}, \tau_\diamond)$.
\end{theorem}

\proof Assume that the sequence $\{x_n\}$   converges to points $x$ and $y$ in $(X, \mathcal{F}, \tau_\diamond)$. Let $\alpha\in (0, 1)$ and $t> 0$. Then there exists $m\in \mathbb{N}$ such that $F_{x_n,\, x}(t)>1-\alpha$ and  $F_{x_n,\, y}(t)>1-\alpha$ for all $n\geq m$. \\
 Now $F_{x,\, y}(t)\geq \tau_\diamond\left( F_{x,\, x_n},\;F_{x_n,\, y}\right)(t)= F_{x,\, x_n}(t)\diamond F_{x_n,\, y}(t) \geq (1-\alpha)\diamond (1-\alpha)\geq (1-\alpha)$ for all $n\geq m$. Therefore $F_{x,\, y}=H_0$ as  $t$ is chosen arbitrarily. So $x=y$.

\begin{theorem}\label{th41}
Let the PM space $(X, \mathcal{F}, \tau_\diamond)$ be  complete and $f: X\rightarrow X$ be a probabilistic TSR-contraction mapping. If $\alpha\diamond\alpha\geq\alpha$ for all $\alpha\in[0, 1]$ then $f$ has a unique fixed point in $X$.\end{theorem}

\proof Since $f$ is a probabilistic TSR-contraction mapping, for $t> 0$ there exists $k\in (0, 1)$ such that $1-F_{f(a), f(b)} (kt)\leq k(1-F_{a, b}(t))$ for all $a, b\in X$.\\
   Let $x_0\in X$ and $x_n=f^n(x_0)$ for all $n\in\mathbb{N}$. Then taking $a=x_1=f(x_0)$ and $b=x_0$ we get
\begin{align*}
    &&1-F_{f(f(x_0)), f(x_0)} (kt)&\leq k(1-F_{f(x_0), x_0}(t)),\\ 
&\text{or}, &1-F_{x_2, x_1} (kt)&\leq k\left(1-F_{x_1, x_0}(t)\right) ,\\
 &\text{or}, &F_{x_2, x_1} (t)&\geq F_{x_2, x_1} (kt)\geq 1-k\left(1-F_{x_1, x_0}(t)\right).\\
&\text{or}, &-F_{x_2, x_1} (t)&\leq -1+k\left(1-F_{x_1, x_0}(t)\right).\\
&\text{or}, &1-F_{x_2, x_1} (t)&\leq k\left(1-F_{x_1, x_0}(t)\right).\\
&\text{or}, &k\left( 1-F_{x_2, x_1} (t)\right) &\leq k^2\left(1-F_{x_1, x_0}(t)\right).\\
&\text{or}, &1-k\left( 1-F_{x_2, x_1} (t)\right) &\geq 1-k^2\left(1-F_{x_1, x_0}(t)\right). \smallskip\hspace{2cm}\cdots (i)
\end{align*}
 Again taking $a=x_2$ and $b=x_1$ we get
 \begin{align*}
   &&1-F_{x_3, x_2} (kt)&\leq k\left( 1-F_{x_2, x_1}(t)\right)  ,\\
  &\text{or}, &F_{x_3, x_2} (t)&\geq F_{x_3, x_2} (kt)\geq 1-k\left( 1-F_{x_2, x_1}(t)\right) \geq 1-k^2\left(1-F_{x_1, x_0}(t)\right).\\
&\text{or}, &-F_{x_3, x_2} (t)&\leq -1+k^2\left(1-F_{x_1, x_0}(t)\right).\\
&\text{or}, &1-F_{x_3, x_2} (t)&\leq k^2\left(1-F_{x_1, x_0}(t)\right).\\
&\text{or}, &k\left( 1-F_{x_3, x_2} (t)\right) &\leq k^3\left(1-F_{x_1, x_0}(t)\right).\\
&\text{or}, &1-k\left( 1-F_{x_3, x_2} (t)\right) &\geq 1-k^3\left(1-F_{x_1, x_0}(t)\right).\smallskip\hspace{2cm}\cdots (ii)
 \end{align*}
Again taking $a=x_3$ and $b=x_2$ we get $1-F_{x_4, x_3} (kt)\leq k\left( 1-F_{x_3, x_2}(t)\right)  $, that is, $F_{x_4, x_3} (t)\geq F_{x_4, x_3} (kt)\geq 1-k\left( 1-F_{x_3, x_2}(t)\right) \geq 1-k^3\left(1-F_{x_1, x_0}(t)\right) $, by $(ii)$. 
So, by induction we get $F_{x_{n+1},\, x_n} (t)\geq 1- k^n\left( 1-F_{x_1,\, x_0}(t)\right) $ for all $n\in\mathbb{N}$.\\
We now show that $\{x_n\}$ is a Cauchy sequence in $(X, \mathcal{F}, \tau_\diamond)$. Let $\alpha\in (0, 1)$. Then there exists a natural number $m$ such that $ k^{m}\left( 1-F_{x_1,\, x_0}(t)\right)< \alpha$. Let $p $ be any natural number. Now,
\begin{align*}
F_{x_{n+p}, \,x_n}(t)&\geq \tau_\diamond\left( F_{x_{n+p}, \,x_{n+p-1}},F_{x_{n+p-1}, \,x_{n}}\right)(t)\\
&= F_{x_{n+p}, \,x_{n+p-1}}(t)\diamond F_{x_{n+p-1}, \,x_{n}}(t)\\
&\geq  F_{x_{n+p}, \,x_{n+p-1}}(t)\diamond \tau_\diamond\left( F_{x_{n+p-1}, \,x_{n+p-2}}, F_{x_{n+p-2}, \,x_{n}}\right)(t), \text{ as } \diamond \text{ is non-decreasing.}\\
&= F_{x_{n+p}, \,x_{n+p-1}}(t)\diamond  F_{x_{n+p-1}, \,x_{n+p-2}}(t)\diamond F_{x_{n+p-2}, \,x_{n}}(t), \text{ as } \diamond \text{ is associative.}\\
&\geq\\
&\vdots\\
&\geq\\
&= F_{x_{n+p}, \,x_{n+p-1}}(t)\diamond  F_{x_{n+p-1}, \,x_{n+p-2}}(t)\diamond \cdots \diamond F_{x_{n+1}, \,x_{n}}(t)\\
&\geq \left( 1- k^{n+p-1}\left( 1-F_{x_1,\, x_0}(t)\right)\right)  \diamond \left( 1- k^{n+p-2}\left( 1-F_{x_1,\, x_0}(t)\right) \right) \diamond \cdots \diamond \\
 &\smallskip\hspace{.5 cm}\left( 1- k^{n}\left( 1-F_{x_1,\, x_0}(t)\right)\right)\\
&\geq \left( 1- k^{n}\left( 1-F_{x_1,\, x_0}(t)\right)\right)  \diamond \left( 1- k^{n}\left( 1-F_{x_1,\, x_0}(t)\right) \right) \diamond \cdots \diamond \\ &\smallskip\hspace{.5 cm}\left( 1- k^{n}\left( 1-F_{x_1,\, x_0}(t)\right)\right), \text{ as } k\in (0, 1) \text{ and } \diamond \text{ is non-decreasing.}\\
&\geq 1- k^{n}\left( 1-F_{x_1,\, x_0}(t)\right).
\end{align*}
So, $F_{x_{n+p}, \,x_n}(t)>1-\alpha$ whenever $1- k^{n}\left( 1-F_{x_1,\, x_0}(t)\right)> 1-\alpha$.\\
Now $1- k^{n}\left( 1-F_{x_1,\, x_0}(t)\right)> 1-\alpha$ iff $ k^{n}\left( 1-F_{x_1,\, x_0}(t)\right)< \alpha$ if $n\geq m$, as\\
 $ k^{m}\left( 1-F_{x_1,\, x_0}(t)\right)< \alpha$.
Thus $\{x_n\}$ is a cauchy sequence in $(X, \mathcal{F}, \tau_\diamond)$. Since the space is complete, there exists $y\in X$ such that $\lim x_n=y$, that is, $\lim f^n (x_0)=y$. We now show that $y$ is a fixed point of the function $f$, that is, $f(y)=y$.  We have
\begin{align*}
F_{f(y),\, y}(t)&\geq \tau_\diamond\left( F_{f(y),\, x_n},\, F_{x_n, y}\right) (t)\\
 &= F_{f(y),\, x_n}(t)\diamond F_{x_n, y}(t)\\
 &=F_{f(y),\, f^n(x_0)}(t)\diamond F_{x_n, y}(t)\\
& = F_{f(y),\, f(x_{n-1})}(t)\diamond F_{x_n, y}(t)\\
& \geq \left[ 1-k\left( 1-F_{y,\, x_{n-1}}\left( \frac{t}{k}\right) \right) \right] \diamond F_{x_n, y}(t).\smallskip\hspace{2cm}\cdots (iii)
 \end{align*}
 Now since $\lim x_n = y$, for any $\alpha \in (0, 1)$, there exists $m\in \mathbb{N}$ such that $F_{y,\, x_{n-1}}(\frac{t}{k})> 1-\alpha$ and $F_{x_n, y}(t)> 1-\alpha$ for all $n\geq m$. So, from $(iii)$ we get \\
 $F_{f(y),\, y}(t)\geq ( 1-k\alpha ) \diamond (1-\alpha)\geq  1-\alpha$ as $\diamond$ is non-decreasing.\\
 So,  $F_{f(y),\, y}(t)\geq  1-\alpha$ for all $t> 0$.\\
 Therefore $F_{f(y),\, y}(t)=1$ for $t>0$. So $F_{f(y),\, y}= H_0$.\\
Thus $f(y)=y$. We now verify that there exists only one fixed point of $f$. Let $z\in X$ be such that $f(z)=z$. Then $F_{y, z}(t)\geq F_{y, z}(kt)= F_{f(y), f(z)}(kt)\geq 1-k\left( 1-F_{y, z} (t)\right)$ implies $1- F_{y, z}(t)\leq k\left( 1-F_{y, z} (t)\right)$. If  $F_{y, z}(t)< 1$, that is, if $1- F_{y, z}(t)> 0$, we obtain $k\geq 1$ which is not true. Thus $F_{y, z}(t)=1$. Since $t$ is chosen arbitrarily, $F_{y, z}= H_0$, that is, $y=z$. This completes the proof of the theorem.

\begin{theorem}\label{th42}
Let us consider the PM space $(X, \mathcal{F}, \tau_\diamond)$ and $\alpha\diamond\alpha\geq\alpha$ for all $\alpha\in[0, 1]$. Let $x_0\in X$, $t> 0$ and $r\in (0, 1)$. Then the t-closed sphere $\overline{S_t}(x_0,\, r,\,\mathcal{F} )$ is t-closed. 
\end{theorem}

\proof Let $y\notin \overline{S_t}(x_0,\, r,\,\mathcal{F} )$. Then $F_{y, x_0}(t)< 1-r$, $\smallskip\hspace{2cm}\cdots (i)$\\
We now show that $S_t(y,\, r,\,\mathcal{F} )\cap \overline{S_t}(x_0,\, r,\,\mathcal{F} )=\phi$. If possible, let $z\in S_t(y,\, r,\,\mathcal{F} )\cap \overline{S_t}(x_0,\, r,\,\mathcal{F} )$. Then $F_{y, z}(t)> 1-r$ and $F_{z, x_0}(t)\geq 1-r$. Now\\
 $F_{y, x_0}(t)\geq \tau_\diamond\left( F_{y, z},\, F_{z, x_0}\right) (t)= F_{y, z}(t)\diamond F_{z, x_0}(t)\geq (1-r)\diamond (1-r)\geq 1-r $, which contradicts $(i)$. Thus $y$ is not a limit point of $\overline{S_t}(x_0,\, r,\,\mathcal{F} )$. Hence $\overline{S_t}(x_0,\, r,\,\mathcal{F} )$ is t-closed.

\begin{theorem}\label{th43}
Let the PM space $(X, \mathcal{F}, \tau_\diamond)$ be complete and $f: X\rightarrow X$ be a probabilistic TSR-contraction mapping on $\overline{S_t}(x_0,\, r,\,\mathcal{F} )$ for some $t>0$. If $\alpha\diamond\alpha\geq\alpha$ for all $\alpha\in[0, 1]$ and moreover $F_{x_0,\, f(x_0)}(u)>1-r$ for all $u> 0$, then $f$ has a unique fixed point in $\overline{S_t}(x_0,\, r,\,\mathcal{F})$.
\end{theorem}

\proof Let $x_1=f(x_0)$, $x_2=f(x_1)=f^2(x_0)$, $\cdots$,  $x_n=f(x_{n-1})=f^n(x_0)$. \\
Now $F_{x_1,\, x_0}(t)=F_{f(x_0),\, x_0}(t)>1-r$. So, $x_1\in \overline{S_t}(x_0,\, r,\,\mathcal{F} )$. Assume that \\
$x_1, x_2, \cdots, X_{n-1}\in \overline{S_t}(x_0,\, r,\,\mathcal{F} )$. We now show that $x_n\in \overline{S_t}(x_0,\, r,\,\mathcal{F} )$.
Since $f$ is a probabilistic TSR contraction mapping on $\overline{S_t}(x_0,\, r,\,\mathcal{F} )$, there exists $k\in (0, 1)$ such that $1-F_{f(x),\, f(y)}(kt)\leq k\left( 1-F_{x,\, y}(t)\right)$ for all $x,\, y\in \overline{S_t}(x_0,\, r,\,\mathcal{F} )$.\\ 
 Since $x_0,\, x_1\in \overline{S_t}(x_0,\, r,\,\mathcal{F} )$,
\begin{align*}
&&1-F_{x_2,\, x_1}(kt)&=1-F_{f(x_1),\, f(x_0)}(kt)\leq k\left( 1-F_{x_1,\, x_0}(t)\right) ,\\
&\text{or}, &1- \frac{1-F_{x_2,\, x_1}(kt)}{k}&\geq F_{x_1,\, x_0}(t)>1-r,  \text{ as } x_1\in \overline{S}(x_0,\, r,\,\mathcal{F} ),\\
&\text{or}, &\frac{1-F_{x_2,\, x_1}(kt)}{k}&<r,\\
&\text{or}, &F_{x_2,\, x_1}(kt)&>1-kr,\\
&\text{or}, &F_{x_2,\, x_1}(t)&>1-kr \smallskip\hspace{2cm}\cdots (i)
\end{align*}
Again, \begin{align*}
&&1-F_{x_3,\, x_2}(kt)&=1-F_{f(x_2),\, f(x_1)}(kt)\leq k\left( 1-F_{x_2,\, x_1}(t)\right),\\
&\text{or}, &1- \frac{1-F_{x_3,\, x_2}(kt)}{k}&\geq F_{x_2,\, x_1}(t)>1-kr,  \text{by} (i),\\
&\text{or}, &\frac{1-F_{x_3,\, x_2}(kt)}{k}&<kr,\\
&\text{or}, &F_{x_3,\, x_2}(kt)&>1-k^2r.\\
&\text{or}, &F_{x_3,\, x_2}(t)&>1-k^2r.
\end{align*}
So, by induction, $F_{x_n,\, x_{n-1}}(t)>1-k^{n-1}r$ for all $n\in \mathbb{N}$.
Now,
\begin{align*}
F_{x_{n}, \,x_0}(t)&\geq \tau_\diamond\left( F_{x_{n}, \,x_{n-1}},F_{x_{n-1}, \,x_{0}}\right)(t)\\
&= F_{x_{n}, \,x_{n-1}}(t)\diamond F_{x_{n-1}, \,x_{0}}(t)\\
&\geq  F_{x_{n}, \,x_{n-1}}(t)\diamond \tau_\diamond\left( F_{x_{n-1}, \,x_{n-2}}, F_{x_{n-2}, \,x_{0}}\right)(t), \text{ as } \diamond \text{ is non-decreasing.}\\
&= F_{x_{n}, \,x_{n-1}}(t)\diamond  F_{x_{n-1}, \,x_{n-2}}(t)\diamond F_{x_{n-2}, \,x_{0}}(t), \text{ as } \diamond \text{ is associative.}\\
&\geq\\
&\vdots\\
&\geq\\
&= F_{x_{n}, \,x_{n-1}}(t)\diamond  F_{x_{n-1}, \,x_{n-2}}(t)\diamond \cdots \diamond F_{x_{1}, \,x_{0}}(t)\\
&\geq \left( 1- k^{n-1}r\right)   \diamond \left( 1- k^{n-2}r \right) \diamond \cdots \diamond \left( 1- kr\right) \diamond \left( 1- r\right)\\
&\geq \left( 1- r\right)   \diamond \left( 1- r \right) \diamond \cdots \diamond \left( 1- r\right) \diamond \left( 1- r\right), \text{ as } k\in (0, 1) \text{ and }\\
&\smallskip\hspace{7cm} \diamond \text{ is non-decreasing.}\\
&\geq 1- r.
\end{align*}
Therefore, $x_n\in \overline{S_t}(x_0,\, r,\,\mathcal{F} )$. Thus, $x_n\in \overline{S_t}(x_0,\, r,\,\mathcal{F} )$ for all $n\in \mathbb{N}$.
Let $u>0$. Since $f$ is a probabilistic TSR contraction mapping on $\overline{S_t}(x_0,\, r,\,\mathcal{F} )$, there exists $p\in (0, 1)$ such that $1-F_{f(x),\, f(y)}(pu)\leq p\left( 1-F_{x,\, y}(u)\right)$ for all $x,\, y\in \overline{S_t}(x_0,\, r,\,\mathcal{F} )$.\\ 
Then by the above procedure, we get  $F_{x_n,\, x_{n-1}}(u)>1-p^{n-1}r$ for all $n\in \mathbb{N}.$\\
Now by the procedure of the proof of the Theorem \ref{th41}, we can show that the sequence $\{x_n\}$ is a Cauchy sequence. Since $(X, \mathcal{F}, \tau_\diamond)$ is a complete PM space,  $\{x_n\}$ converges to an element $y$ of $X$. So, $y$ is a t-limit point of $\overline{S_t}(x_0,\, r,\,\mathcal{F} )$. Again Since  $\overline{S_t}(x_0,\, r,\,\mathcal{F} )$ is t-closed, $y\in \overline{S_t}(x_0,\, r,\,\mathcal{F} )$.  We now show that $y$ is a unique fixed point of $f$ in $\overline{S_t}(x_0,\, r,\,\mathcal{F} )$. Let $z$ be a fixed point of $f$ in $\overline{S_t}(x_0,\, r,\,\mathcal{F} )$. Then $F_{y, z}(t)\geq F_{y, z}(pu)= F_{f(y), f(z)}(pu)\geq 1-p\left( 1-F_{y, z} (u)\right)$ implies $1- F_{y, z}(u)\leq p\left( 1-F_{y, z} (u)\right)$. If  $F_{y, z}(u)< 1$, that is, if $1- F_{y, z}(u)> 0$, we obtain $p\geq 1$ which is not true. Thus $F_{y, z}(u)=1$. Since $u$ is chosen arbitrarily, $F_{y, z}= H_0$, that is, $y=z$. This completes the proof of the theorem. 

\begin{theorem}\label{th44}
Let us consider the PM space $(X, \mathcal{F}, \tau_\diamond)$ and $\alpha\diamond\alpha\geq\alpha$ for all $\alpha\in[0, 1]$. If $x_n\rightarrow x$ and $y_n\rightarrow y$ as $n\rightarrow\infty$ in $(X, \mathcal{F}, \tau_\diamond)$ then $F_{x_n,\, y_n}(t)\rightarrow F_{x,\, y}(t) $ for all $t> 0$.
\end{theorem}

\proof Let $t> 0$ and $F_{x, y}(t)=\alpha$. Since $x_n\rightarrow x$ and $y_n\rightarrow y$, for any $\epsilon\in (0, \frac{\alpha}{2})$, there exists $m\in\mathbb{N}$ such that $F_{x_n,\, x}(t)>1-\epsilon$ and $F_{y_n,\, y}(t)>1-\epsilon$ for all $n\geq m$. Now 
\begin{align*}
F_{x_n,\, y_n}(t)&\geq \tau_\diamond\left( F_{x_n,\, x}, F_{x, y_n}\right) (t)\\
&= F_{x_n,\, x}(t)\diamond F_{x, y_n}(t)\\
&\geq F_{x_n,\, x}(t)\diamond \tau_\diamond\left( F_{x, y}, F_{y, y_n}\right) (t)\\
&=F_{x_n,\, x}(t)\diamond F_{x, y}(t)\diamond F_{y, y_n}(t)\\
&\geq (1-\epsilon)\diamond F_{x, y}(t)\diamond (1-\epsilon) \text{ for all } n\geq m.\\
&\geq (1-\epsilon)\diamond F_{x, y}(t)\\
&\geq (1-\epsilon)\diamond \left( F_{x, y}(t)-\epsilon\right), \text{ as } F_{x, y}(t)>F_{x, y}(t)-\epsilon=\alpha-\epsilon> 0\\
&\geq \left( F_{x, y}(t)-\epsilon\right)\diamond \left( F_{x, y}(t)-\epsilon\right), \text{ as } 1\geq F_{x, y}(t)\\
&\geq F_{x, y}(t)-\epsilon.
\end{align*}
Therefore $F_{x_n,\, y_n}(t)-F_{x, y}(t) \geq -\epsilon$ for all $n\geq m$ and for all $t\in\mathbb{R}$.
Again, 
\begin{align*}
F_{x,\, y}(t)&\geq \tau_\diamond\left( F_{x,\, x_n}, F_{x_n, y}\right) (t)\\
&= F_{x,\, x_n}(t)\diamond F_{x_n, y}(t)\\
&\geq F_{x_n,\, x}(t)\diamond \tau_\diamond\left( F_{x_n, y_n}, F_{y_n, y}\right) (t)\\
&=F_{x_n,\, x}(t)\diamond F_{x_n, y_n}(t)\diamond F_{y_n, y}(t)\\
&\geq (1-\epsilon)\diamond F_{x_n, y_n}(t)\diamond (1-\epsilon) \text{ for all } n\geq m.\\
&\geq (1-\epsilon)\diamond F_{x_n, y_n}(t)\\
&\geq (1-\epsilon)\diamond \left( F_{x_n, y_n}(t)-\epsilon\right),\\
&\smallskip\hspace{1cm} \text{ as } F_{x_n, y_n}(t)\geq \left( F_{x_n, y_n}(t)-\epsilon\right) \geq F_{x, y}(t)-2\epsilon=\alpha -2\epsilon> 0,\\
&\geq \left( F_{x_n, y_n}(t)-\epsilon\right)\diamond \left( F_{x_n, y_n}(t)-\epsilon\right), \text{ as } 1-\epsilon\geq F_{x_n, y_n}(t)-\epsilon,\\
&\geq  F_{x_n, y_n}(t)-\epsilon.
\end{align*}
Therefore $F_{x_n, y_n}(t)-F_{x, y}(t)\leq \epsilon$.
Thus, $-\epsilon \leq F_{x_n, y_n}(t)-F_{x, y}(t)\leq \epsilon$, that is, $|F_{x_n, y_n}(t)-F_{x, y}(t)|\leq\epsilon$ for all $t\in \mathbb{R}.$ Hence $F_{x_n, y_n}\rightarrow F_{x, y}$ as $n\rightarrow\infty$

\begin{theorem}
Let the PM space $(X, \mathcal{F}, \tau_\diamond)$ be complete, $f: X\rightarrow X$ be a continuous function and $\alpha\diamond\alpha\geq\alpha$ for all $\alpha\in[0, 1]$. If $f^m$ is a probabilistic TSR-contraction mapping for some natural number $m$ then $f$ has a unique fixed point in $X$.
\end{theorem}

\proof Let $g=f^m$. Since $f^m$ is a contraction mapping, for $t> 0$ there exists $k\in (0, 1)$ such that $1-F_{g(x),\; g(y)}(kt)\leq k(1-F_{x,\, y})(t)$  that is, $F_{g(x),\; g(y)}(kt)\geq 1- k(1-F_{x,\, y}(t))$.\\
If $n$ is a positive integer and $x_0\in X$ then
\begin{align*}
F_{g^nf(x_0),\;g^n(x_0)}(t)&\geq F_{g^nf(x_0),\;g^n(x_0)}(kt)\\
&=F_{g\left( g^{n-1}f(x_0)\right),\;g\left( g^{n-1}(x_0)\right)} (kt)\\
&\geq  1-k\left[ 1-F_{ g^{n-1}f(x_0),\; g^{n-1}(x_0)} (t)\right]\\
&= 1-k\left[ 1-F_{ gg^{n-2}f(x_0),\; gg^{n-2}(x_0)} (t)\right].
\end{align*}
Now 
\begin{align*}
&&F_{ gg^{n-2}f(x_0),\; gg^{n-2}(x_0)} (t)&\geq F_{ gg^{n-2}f(x_0),\; gg^{n-2}(x_0)} (kt)\\
&&&\geq 1-k\left[ 1-F_{ g^{n-2}f(x_0),\; g^{n-2}(x_0)} (t)\right]\\ 
&\text{or}, &1- F_{ gg^{n-2}f(x_0),\; gg^{n-2}(x_0)} (t)&\leq k\left[ 1-F_{ g^{n-2}f(x_0),\; g^{n-2}(x_0)} (t)\right]\\
&\text{or}, &1-k\left[ 1-F_{ gg^{n-2}f(x_0),\; gg^{n-2}(x_0)} (t)\right]&\geq 1- k^2\left[ 1-F_{ g^{n-2}f(x_0),\; g^{n-2}(x_0)} (t)\right].
\end{align*}
Therefore 
\begin{align*}
F_{g^nf(x_0),\;g^n(x_0)}(t)&\geq 1-k\left[ 1-F_{ gg^{n-2}f(x_0),\; gg^{n-2}(x_0)} (t)\right]\\
&\geq 1- k^2\left[ 1-F_{ g^{n-2}f(x_0),\; g^{n-2}(x_0)} (t)\right]\\
&\geq\\ &\cdots\\ &\geq 1- k^n\left[ 1-F_{ f(x_0),\; x_0} (t)\right].
\end{align*}
Now $ K^n\rightarrow 0$ as $n\rightarrow \infty$.\\
So, $F_{g^nf(x_0),\;g^n(x_0)}(t)\rightarrow 1$ for all $t> 0$ as $n\rightarrow \infty$, $\smallskip\hspace{2cm}\cdots (i)$\\
Now by the Theorem \ref{th41}, $g$ has a unique fixed point, say $y$ and $\lim g^n(x_0)=y$. Since $f$ is continuous, by the Theorem \ref{th21}, $\lim fg^n(x_0)=f(y)$, that is, $\lim g^nf(x_0)=f(y)$.
So, by the Theorem \ref{th44}, $F_{g^nf(x_0),\; g^n(x_0)}(t)\rightarrow F_{f(y),\, y}(t)$ as $n\rightarrow\infty$,\\
 Therefore $F_{g^nf(x_0),\; g^n(x_0)}(t)\rightarrow F_{f(y),\, y}(t)$ for all $t\in \mathbb{R}$ as $n\rightarrow\infty$.\\
 Since every convergent sequence has a unique limit, $F_{f(y),\, y}(t)=1$ for all $t>0$. Therefore $F_{f(y),\, y}= H_0$ which implies that $f(y)=y$. So, $y$ is a fixed point of $f$.\\
 Now if $y_1$ is a fixed point of $f$, that is, $f(y_1)=y_1$. Then $f^m(y_1)=f^{m-1}f(y_1)=f^{m-1}(y_1)=\cdots=y_1$. So, $y_1$ is a fixed point of $f^m=g$. Since $g$ has a unique fixed point $y$, $y_1=y$. This completes the proof of the theorem.

\end{document}